\documentclass{article}
\usepackage[utf8]{inputenc}
\usepackage{amsmath}
\usepackage{amssymb}
\usepackage{amsthm}
\usepackage{enumerate}
\usepackage{mathrsfs}
\usepackage{scrextend}
\usepackage{mathtools}
\usepackage[hyperfootnotes=false]{hyperref}
\usepackage{cleveref}
\usepackage[multiple]{footmisc}
\usepackage[a4paper, total={6in, 9in}]{geometry}
\usepackage{scalerel,stackengine}
\usepackage{tikz}
\usepackage{float}
\stackMath
\newcommand\reallywidehat[1]{%
\savestack{\tmpbox}{\stretchto{%
  \scaleto{%
    \scalerel*[\widthof{\ensuremath{#1}}]{\kern-.6pt\bigwedge\kern-.6pt}%
    {\rule[-\textheight/2]{1ex}{\textheight}}
  }{\textheight}%
}{0.5ex}}%
\stackon[1pt]{#1}{\tmpbox}%
}
\theoremstyle{plain}
\newtheorem{theorem}{Theorem}
\newtheorem{lemma}[theorem]{Lemma}
\newtheorem{proposition}[theorem]{Proposition}
\newtheorem{corollary}[theorem]{Corollary}

\theoremstyle{definition}

\newtheorem{question}[theorem]{Question}

\numberwithin{theorem}{section}

\title{Improved bounds for skew corner-free sets}
\author{Adrian Beker\footnote{University of Zagreb, Faculty of Science, Department of Mathematics, Zagreb,
Croatia.\\ Email: \nolinkurl{adrian.beker@math.hr}}}
\date{\today}

\begin{document}

\maketitle

\begin{abstract}
    We construct skew corner-free subsets of $[n]^2$ of size $n^2\exp(-O(\sqrt{\log n}))$, thereby improving on recent bounds of the form $\Omega(n^{5/4})$ obtained by Pohoata and Zakharov. In the other direction, we prove that any such set has size at most $O(n^2(\log n)^{-c})$ for some absolute constant $c > 0$. This improves on the previously best known upper bound, coming from Shkredov's work on the corners theorem.
\end{abstract}

\section{Introduction}

A \emph{skew corner} is a configuration of points of the form $(x,y), (x, y + d), (x+d,y')$ with $x,y,y',d\in \mathbb{Z}$. We say that a skew corner is \emph{trivial} if $d = 0$ and \emph{non-trivial} otherwise. A subset of $\mathbb{Z}^2$ is said to be \emph{skew corner-free} if it contains no non-trivial skew corners. These notions were recently introduced by Pratt \cite{pratt} (see also \cite{pratt-mathoverflow}), who asked the following interesting question: what is the largest size of a skew corner-free subset of $[n]^2$? 

Letting $s(n)$ denote the answer to this question and taking $y' = y$ in the definition of a skew corner, one immediately sees that $s(n)$ cannot exceed the largest size of a corner-free subset of $[n]^2$. The latter quantity, which we will denote by $r(n)$, is an important object of study in additive combinatorics. Ajtai and Szemerédi \cite{ajtai-szemeredi} were the first to establish that $r(n) = o(n^2)$, a result now known as the \emph{corners theorem}. A lot of work has since been done on this topic, see e.g.\ \cite{solymosi} for an elegant proof using the triangle removal lemma. The best known quantitative bounds on $r(n)$ are due to Shkredov \cite{shkredov} and Green \cite{green-corners}, who showed that $r(n) = O(n^2(\log\log n)^{-1/73})$ and $r(n) \geq 2^{-(c+o(1))\sqrt{\log_2n}}n^2$, where $c = 2\sqrt{2\log_2\frac{4}{3}}$, respectively. For a version of Shkredov's argument in the finite field setting, which preserves the main ideas but greatly simplifies the technical details, we refer the reader to \S5 of Green's excellent survey \cite{green-finite-fields} and the companion note \cite{green-shkredov}.

Returning to skew corners, the additional degree of freedom in the definition makes it seem rather hard to construct large sets lacking these structures. This led Pratt to conjecture that, for any $\varepsilon > 0$, one has $s(n) = O(n^{1+\varepsilon})$ (Conjecture 1.2 in \cite{pratt}). The significance of this conjecture is that, were it to be true, it would rule out the the possibility of obtaining algorithms for matrix multiplication of complexity $O(n^{2+o(1)})$ using certain group-theoretic approaches. However, this prediction turned out to be false -- in very recent work \cite{pohoata-zakharov}, Pohoata and Zakharov showed that $s(n)$ grows at least as fast as $n^{5/4}$. Nonetheless, they hypothesised that one should still have $s(n) = O(n^{2-c})$ for some absolute constant $c > 0$. In this paper, we prove to the contrary that there exist skew corner-free subsets of $[n]^2$ of size $n^{2-o(1)}$. More precisely, we prove the following lower bound, which has the same shape as the corresponding bound for corner-free sets.

\begin{theorem}
\label{main lower bound}
There exists a constant $c_1 > 0$ such that
$$s(n) \geq \frac{n^2}{2^{(c_1+o(1))\sqrt{\log_2n}}}.$$
\end{theorem}

The proof of Theorem \ref{main lower bound} is based on a Behrend-type construction (see \cite{behrend}). In Section 2, we first give an example which suffices to disprove Pratt's conjecture, but is, in our opinion, somewhat simpler than that of \cite{pohoata-zakharov}. We then proceed to give a full proof of Theorem \ref{main lower bound} and briefly comment on how it relates to the arguments of \cite{pohoata-zakharov}.

In the other direction, it is natural to wonder whether one can improve the bound on $s(n)$ obtained from Shkredov's bounds on $r(n)$. Furthermore, Shkredov's argument is rather involved, so it would be desirable to have a simpler way of obtaining bounds for skew corner-free sets. We do so by proving the following result.

\begin{theorem}
\label{main upper bound}
There exists a constant $c_2 > 0$ such that 
$$s(n) = O\Biggl(\frac{n^2}{(\log n)^{c_2}}\Biggr).$$
\end{theorem}

The proof of Theorem \ref{main upper bound} is based on the traditional Fourier-analytic approach, as opposed to \cite{shkredov}, which uses a different, more subtle notion of uniformity. More precisely, we adapt the approach of Heath-Brown \cite{heath-brown} and Szemerédi \cite{szemeredi} to Roth's theorem. For a very readable exposition of these arguments, the reader may consult \cite{green-szemeredi} (see also \cite{peluse}). Even though the general scheme of the argument transfers fairly smoothly to our setting, we have to make several innovations on a technical level in order to make it work. This is done in Section 3, whereas in Section 4, we discuss how this problem relates to other Szemerédi-type problems and highlight some questions that remain open.

It will be useful to consider skew corner-free sets in the setting of abelian groups other than $\mathbb{Z}$, namely $\mathbb{Z}^d$ and finite cyclic groups. An observation that turns out to be convenient in proofs of both lower and upper bounds is that a subset of $[n]^2$ retains the property of being skew corner-free if we regard it as a subset of $(\mathbb{Z}/2n\mathbb{Z})^2$. Another feature of this property, of which we make crucial use in the proof of Theorem \ref{main upper bound}, is that it is translation-invariant in a rather strong sense. Namely, if $G$ is a finite cyclic group then the skew corner-freeness of a subset of $G^2$ is preserved under horizontal and vertical translations, where we are allowed to perform a different translation inside each column of the toroidal grid $G^2$. This fact will enable us to carry out the density increment strategy similarly as in the proof of Roth's theorem.

\bigskip

\noindent\textbf{Notation and normalisations.} We use standard asymptotic notation. Given quantities $A$ and $B$ that depend (in a possibly implicit way) on a positive integer parameter $n$, we write $A = O(B)$ or $A \ll B$ if there is a constant $K > 0$ such that $|A| \leq K|B|$ for all sufficiently large $n$. We also use the equivalent notation $B = \Omega(A)$ and $B \gg A$. We write $A = o(B)$ to mean that $\frac{A}{B} \to 0$ as $n \to \infty$.

We now introduce some notation and conventions pertaining to Fourier analysis. Let $G$ be a finite abelian group. Given a positive integer $k$ and a function $f\colon G^k \to \mathbb{C}$, we denote by
$$\mathbb{E}_{x_1, \ldots, x_k\in G}f(x_1,\ldots,x_k) \vcentcolon= \frac{1}{|G^k|}\sum_{x\in G^k}f(x)$$
its average over $G^k$. Given functions $f,g\colon G \to \mathbb{C}$, we define their inner product by
$$\langle f,g \rangle \vcentcolon= \mathbb{E}_{x\in G}f(x)\overline{g(x)}.$$
For any $p \in [1, \infty)$, we define the $\mathbb{L}^p$-norm by
$$\lVert f\rVert_p \vcentcolon= \Bigl(\mathbb{E}_{x\in G}|f(x)|^p\Bigr)^{\frac{1}{p}},$$
whereas the $\mathbb{L}^{\infty}$ norm is defined in the usual way as
$$\lVert f\rVert_{\infty} \vcentcolon= \max_{x\in G}|f(x)|.$$
We define the difference convolution of $f$ and $g$ to be
$$f \circ g \colon G \to \mathbb{C}, \quad x \mapsto \mathbb{E}_{y\in G}f(y)\overline{g(y-x)}.$$
We write $\widehat{G}$ for the group of characters of $G$ and define the Fourier transform of $f$ by
$$\widehat{f} : \widehat{G} \to \mathbb{C}, \quad \gamma \mapsto \langle f,\gamma\rangle = \mathbb{E}_{x\in G}f(x)\overline{\gamma(x)}.$$
We will make use of the Fourier inversion formula
$$f(x) = \sum_{\gamma \in \widehat{G}}\widehat{f}(\gamma)\gamma(x),$$
Parseval's identity
$$\mathbb{E}_{x\in G}f(x)\overline{g(x)} = \sum_{\gamma \in \widehat{G}}\widehat{f}(\gamma)\overline{\widehat{g}(\gamma)},$$
and the fact that $\reallywidehat{f \circ g} = \widehat{f}\overline{\widehat{g}}$. When working with characters in explicit form, we will find it convenient to use the standard notation $e(\theta) \vcentcolon= e^{2\pi i\theta}$ and $\lVert \theta \rVert_{\mathbb{R}/\mathbb{Z}} \vcentcolon= \min\{|x| \mid x \in \theta\}$ for $\theta \in \mathbb{R}/\mathbb{Z}$, as well as the elementary estimates 
$$4\lVert \theta \rVert_{\mathbb{R}/\mathbb{Z}} \leq |1-e(\theta)| = 2|\sin(\pi\theta)| \leq 2\pi\lVert \theta \rVert_{\mathbb{R}/\mathbb{Z}}.$$
We next introduce some notation specific to our needs. Given a function $f\colon G^2 \to \mathbb{C}$, we define
$$S(f) \colon G \to \mathbb{C}, \quad x \mapsto \mathbb{E}_{y\in G}f(x,y).$$
Thus, $S(f)(x)$ coincides with $\lVert f(x,\cdot)\rVert_1$ if $f$ is non-negative. We also define the following normalised variant of $f$:
$$\widetilde{f} \colon G^2 \to \mathbb{C}, \quad (x,y) \mapsto \begin{cases}\frac{f(x,y)}{\lVert f(x,\cdot)\rVert_1} & \text{if } f(x,\cdot) \neq 0\\0 & \text{otherwise}\end{cases}.$$
Finally, in Section 2, we will be using the notation $\langle \cdot, \cdot \rangle$ and $\lVert \cdot \rVert$ to denote the standard inner product and Euclidean norm on $\mathbb{R}^d$ respectively. These should not be confused with the notational conventions previously laid out in this section since they are defined with respect to the counting measure instead of the uniform probability measure.

\section{Proof of Theorem \ref{main lower bound}}

We first present a simple proof of the weaker bound $s(n) \gg n^{1+c_1}$ for some constant $c_1 > 0$. We will be fairly brief since this is superseded by the results of \cite{pohoata-zakharov}. The idea is to choose a fixed skew corner-free set $S \subseteq (\mathbb{Z}/b\mathbb{Z})^2$ of size $|S| > b$, for some positive integer $b$. That one can do so follows from the discussion at the end of the Section 1 since it is known that $s(n)$ is not $O(n)$ (see Proposition 4.16 in \cite{pratt}). Given a positive integer $n$, choose $k = \lfloor \log_bn\rfloor$ and let $A \subseteq [n]^2$ consist of all points $(x,y) \in [b^k]^2$ such that $(x_j,y_j) \in S$ for all $j \in [k]$, where $z_j$ denotes the $j$-th digit in the base-$b$ representation of $z-1$ for $z \in [b^k]$. It is straightforward to check that $A$ is skew corner-free and that it has size $\Omega(n^{1+c_1})$ with $c_1 = \frac{1}{2}(\log_b|S|-1) > 0$ say. We remark that this kind of product construction is analogous to the binary-ternary Salem-Spencer-type construction for sets free of three-term arithmetic progressions (see \cite{salem-spencer}).

We now turn to the main business of this section, which is to establish Theorem \ref{main lower bound}. Given a positive integer $n$, consider the box $B = [m]^d \subseteq \mathbb{Z}^d$, where $m,d$ are positive integers such that $(2m)^d \leq n$. For positive integers $r,t$, consider the set
$$A_{r,t} = \{(x,y) \in B \times B \mid \lVert x\rVert^2 = r,\ \langle x,y\rangle = t\}.$$
We claim that $A_{r,t}$ is skew corner-free as a subset of $(\mathbb{Z}^d)^2$. Indeed, suppose that it contains three points of the form $(x,y), (x,y''), (x',y')$ with $y'' - y = x' - x$. Since $\langle x,y\rangle = t = \langle x,y''\rangle$, it follows that $\langle x,y''-y\rangle = 0$, whence $\langle x,x'-x\rangle = 0$. Hence, Pythagoras' theorem implies that
$$r = \lVert x'\rVert^2 = \lVert x\rVert^2 + \lVert x'-x\rVert^2 = r + \lVert x'-x\rVert^2.$$
It follows that $\lVert x'-x\rVert = 0$, that is, $x' = x$, as desired. In order to transfer this example to the integer setting, we use the map
$$\varphi \colon \mathbb{Z}^d \to \mathbb{Z}, \quad (x_1, \ldots, x_d) \mapsto 1+\sum_{j=1}^{d}(2m)^{j-1}(x_j-1),$$
which, when restricted to $B$, becomes a Freiman isomorphism into $[n]$. Thus, $(\varphi \times \varphi)(A_{r,t})$ is a skew corner free-subset of $[n]^2$ of the same size as $A_{r,t}$. Therefore, all that remains is to choose $r,t$ so that $A_{r,t}$ is of the required size. But as $r,t$ range over $[dm^2]$, the sets $A_{r,t}$ form a partition of $B \times B$. Hence, by the pigeonhole principle, there exists a choice of $r,t$ such that
$$|A_{r,t}| \geq \frac{|B\times B|}{(dm^2)^2} = \frac{m^{2d}}{d^2m^4} = \frac{m^{2d-4}}{d^2}.$$
Finally, choosing $m = \lfloor \frac{1}{2}n^{1/d}\rfloor$ and $d = \lfloor (2\log_2n)^{1/2}\rfloor$, it is not hard to check that we get 
$$|A_{r,t}| \geq \frac{n^2}{2^{(4\sqrt{2}+o(1))\sqrt{\log_2n}}}.$$
\textbf{Remark.} Our construction turns out to be quite similar in spirit to that of \cite{pohoata-zakharov}. The construction in \cite{pohoata-zakharov} is based on (the affine version of) the Hermitian unital over $\mathbb{F}_{p^2}^2$, which can be thought of as an analogue of a sphere in the finite field setting. By working directly in the integer setting, we take the advantage of being able to consider higher dimensional spheres, which results in superior bounds.

\noindent\textbf{Remark.} Both constructions described in this section can be easily modified so as to produce bi-skew corner-free sets of comparable size (for a precise definition, see Definition 4.17 in \cite{pratt}). Indeed, for the first construction, one merely has to replace the set $S$ with a suitable bi-skew corner-free set, e.g.\
$$\{(0,0), (0,1), (2, 0), (2, 3), (3, 1), (3, 3), (3, 5), (4, 0)\} \subseteq (\mathbb{Z}/6\mathbb{Z})^2.$$
For the second construction, one can modify the definition of the set $A_{r,t}$ by requiring that the second point also lie on the sphere of radius $\sqrt{r}$ centred at the origin:
$$A_{r,t} = \{(x,y) \in B \times B \mid \lVert x\rVert^2 = \lVert y\rVert^2 = r,\ \langle x,y\rangle = t\}.$$
To ensure that $A_{r,t}$ is large, one first uses the pigeonhole principle to find a value of $r$ such that $\sqrt{r}S^{d-1}$ contains at least $\frac{m^d}{dm^2} = \frac{m^{d-2}}{d}$ points of $B$. A further application of the pigeonhole principle then yields a set $A_{r,t}$ of size at least
$$\frac{(m^{d-2}/d)^2}{dm^2} = \frac{m^{2d-6}}{d^3},$$
which can then be optimised to give a bi-skew corner-free subset of $[n]^2$ of size at least $\frac{n^2}{2^{(4\sqrt{3}+o(1))\sqrt{\log_2n}}}$.

\section{Proof of Theorem \ref{main upper bound}}

Let $G$ be a finite abelian group. We begin by defining the trilinear form
$$\Lambda(f, g, h) = \mathbb{E}_{x,y,y',d\in G}f(x,y)g(x,y+d)h(x+d,y'),$$
where $f,g,h$ are real-valued functions on $G^2$. In particular, if $A \subseteq G^2$, then $N^4\Lambda(1_A,1_A,1_A)$ is the number of skew corners in $A$. By employing the Fourier inversion formula and the fact that the Fourier transform diagonalises convolution, it is straightforward to derive the Fourier representation of $\Lambda$:
\begin{align*}
    \Lambda(f,g,h) &= \mathbb{E}_{x,d\in G}(g(x,\cdot) \circ f(x,\cdot))(d)S(h)(x+d)\\
    &= \mathbb{E}_{x\in G}[S(h)\circ(g(x,\cdot)\circ f(x,\cdot))](x)\\
    &= \mathbb{E}_{x\in G}\sum_{\gamma\in\widehat{G}}\reallywidehat{S(h)\circ(g(x,\cdot)\circ f(x,\cdot))}(\gamma)\gamma(x)\\
    &= \sum_{\gamma\in\widehat{G}}\reallywidehat{S(h)}(\gamma)\mathbb{E}_{x\in G}\reallywidehat{f(x,\cdot)}(\gamma)\overline{\reallywidehat{g(x,\cdot)}(\gamma)}\gamma(x).
\end{align*}
In particular, the triangle inequality implies that
\begin{equation}\label{fourier representation}
    |\Lambda(f,g,h)| \leq \sum_{\gamma\in\widehat{G}}|\reallywidehat{S(h)}(\gamma)|\mathbb{E}_{x\in G}|\reallywidehat{f(x,\cdot)}(\gamma)||\reallywidehat{g(x,\cdot)}(\gamma)|.
\end{equation}

To motivate the argument, we first briefly discuss some of the core ideas in the finite field setting. The key observation is that we can use Fourier analysis to control the count of skew corners in a given set. This is made precise by the following result, which we call a `generalised von Neumann theorem', borrowing the terminology from \cite{green-finite-fields}. Given a parameter $\eta \in [0,1]$, we say that a subset of a finite abelian group is \emph{$\eta$-uniform} if all non-trivial Fourier coefficients of its indicator function are at most $\eta$ in absolute value.

\begin{proposition}[Generalised von Neumann theorem]
\label{generalised von neumann}
Let $G$ be a finite abelian group and let $A \subseteq G^2$ be a subset of density $\alpha$. Let $f_A = 1_A - \alpha1_{G^2}$ be the balanced function of $A$. Then
\begin{equation}\label{controlling the count of skew corners}
    \Lambda(1_A,1_A,1_A) \geq \alpha^3 - \alpha\lVert \reallywidehat{S(f_A)}\rVert_{\infty}.
\end{equation}
In particular, if $A$ is $\eta$-uniform, then $A$ contains at least $(\alpha^3-\alpha\eta)N^4$ skew corners.
\end{proposition}
\noindent\textit{Proof.} Consider the decomposition
$$\Lambda(1_A,1_A,1_A) = \Lambda(1_A,1_A,f_A) + \alpha\Lambda(1_A,1_A,1_{G^2}).$$
By the Cauchy-Schwarz inequality, we have the following lower bound for the second term on the right-hand side:
\begin{equation*}
\Lambda(1_A,1_A,1_{G^2}) = \mathbb{E}_{x,y,d\in G}1_A(x,y)1_A(x,y+d) = \mathbb{E}_{x\in G}S(1_A)(x)^2 \geq \Bigl(\mathbb{E}_{x\in G}S(1_A)(x)\Bigr)^2 = \alpha^2.    
\end{equation*}
On the other hand, by (\ref{fourier representation}) and Parseval's identity, the first term can be upper bounded as follows:
\begin{align*}
    |\Lambda(1_A,1_A,f_A)| &\leq \sum_{\gamma\in\widehat{G}}|\reallywidehat{S(f_A)}(\gamma)|\mathbb{E}_{x\in G}|\reallywidehat{1_A(x,\cdot)}(\gamma)|^2\\
    &\leq \lVert \reallywidehat{S(f_A)} \rVert_{\infty}\mathbb{E}_{x\in G}\sum_{\gamma\in\widehat{G}}|\reallywidehat{1_A(x,\cdot)}(\gamma)|^2\\
    &= \lVert \reallywidehat{S(f_A)} \rVert_{\infty}\mathbb{E}_{x\in G}\mathbb{E}_{y\in G}|1_A(x,y)|^2\\
    &= \alpha\lVert\reallywidehat{S(f_A)}\rVert_{\infty},
\end{align*}
so (\ref{controlling the count of skew corners}) follows. The second statement follows on observing that $\reallywidehat{S(f_A)}(\gamma) = \reallywidehat{f_A}(\gamma,1)$ for $\gamma \in \widehat{G}$, where we naturally identify $\reallywidehat{G^2}$ with $\reallywidehat{G}^2$. $\qed$

Using Proposition \ref{generalised von neumann} as a starting point, one can give a proof of the result in the finite field setting in the same vein as Meshulam's theorem \cite{meshulam}. Indeed, if $G = \mathbb{F}_2^n$ and $A \subseteq G^2$ is skew corner-free of density $\alpha$, then one can easily verify that $\Lambda(1_A,1_A,1_A) \leq \frac{\alpha}{N}$, where $N \vcentcolon= 2^n$. Thus, either $\alpha \leq \sqrt{\frac{2}{N}}$ or there exists $\gamma \in \widehat{G}\setminus\{1\}$ such that $|\reallywidehat{S(f_A)}(\gamma)| \geq \frac{\alpha^2}{2}$. The latter case leads to a density increment of $S(1_A)$ on a coset of a codimension-$1$ subspace $H \leqslant G$, and an extra averaging argument produces a density increment of $A$ on a coset of $H^2$. This density increment can then be iterated to obtain a bound of the form $|A| \ll \frac{4^n}{n}$. This is already better than the best known bound \cite{lacey-mcclain} for corner-free sets in this setting, and the argument is substantially simpler.

However, we do not have the luxury of working in the finite field setting, so we have to work a bit harder. To begin, we fix some notation that will be used throughout the remainder of this section. Let $A \subseteq [n]^2$ be a subset of density $\alpha$. Set $N = 2n$ and view $[n]^2$ as a subset of the group $G^2$, where $G = \mathbb{Z}/N\mathbb{Z}$. Let $f_A = 1_A - \alpha1_{[n]^2}$ be the balanced function of $A$. Given $x \in G$, we let $A_x \vcentcolon= \{y \in G \mid (x,y) \in A\}$ be the slice of $A$ at $x \in G$. We also introduce a large global constant $C > 0$, to be determined later.

In the described setting, the following result plays a role analogous to that of Proposition \ref{generalised von neumann}. Roughly speaking, it states that a skew corner-free set must either have small density or large Fourier coefficients. The key difference is that we do not exploit this fact by deducing the existence of a single large Fourier coefficient, but rather a large $\mathbb{L}^2$-mass captured by several Fourier coefficients. In doing so, for technical reasons that will become clear later, we have to replace one copy of $1_A$ by the normalised variant $\widetilde{1_A}$.

\begin{proposition}
\label{main dichotomy}
Let $A \subseteq [n]^2$ be skew corner-free of density $\alpha$. We have the following dichotomy:
\begin{enumerate}[(i)]
    \item either $\alpha \leq \frac{8}{n}$;
    \item or $\sum_{\gamma \in \widehat{G}\setminus\{1\}}|\reallywidehat{S(f_A)}(\gamma)|\mathbb{E}_{x\in G}|\reallywidehat{\widetilde{1_A}(x,\cdot)}(\gamma)||\reallywidehat{1_A(x,\cdot)}(\gamma)| \geq \frac{\alpha^2}{64}$.
\end{enumerate}
\end{proposition}
\noindent\textit{Proof.} As in the proof of Proposition \ref{generalised von neumann}, we decompose
\begin{equation}\label{main decomposition}
    \Lambda(\widetilde{1_A}, 1_A, 1_A) = \Lambda(\widetilde{1_A}, 1_A, f_A) + \alpha\Lambda(\widetilde{1_A}, 1_A, 1_{[n]^2}).
\end{equation}
Since $A$ is skew corner-free when viewed as a subset of $G^2$, the left-hand side equals
\allowdisplaybreaks
\begin{align*}
    \mathbb{E}_{x,y,y',d\in G}\widetilde{1_A}(x,y)1_A(x,y+d)1_A(x+d,y') &= \frac{1}{N}\mathbb{E}_{x,y,y'\in G}\widetilde{1_A}(x,y)1_A(x,y)1_A(x,y')\\
    &= \frac{1}{N}\mathbb{E}_{x\in G}\Bigl(\mathbb{E}_{y\in G}\widetilde{1_A}(x,y)\Bigr)\Bigl(\mathbb{E}_{y'\in G}1_A(x,y')\Bigr)\\
    &= \frac{1}{N}\mathbb{E}_{x\in G}S(1_A)(x)\\
    &= \frac{\alpha}{4N}.
\end{align*}
To obtain a lower bound on the second term on the right-hand side, we write it as follows:
\begin{align*}
    \alpha\mathbb{E}_{x,y,y',d\in G}\widetilde{1_A}(x,y)1_A(x,y+d)1_{[n]^2}(x+d,y') &= \frac{\alpha}{2}\mathbb{E}_{x,y,d\in G}\widetilde{1_A}(x,y)1_A(x,y+d)1_{[n]}(x+d)\\
    &= \frac{\alpha}{2N^2}\sum_{x\in [n],\ A_x \neq \varnothing}\frac{|\{(y_1,y_2) \in A_x^2 \mid y_2-y_1 \in [n]-x \}|}{|A_x|}.
\end{align*}
Since all pairs $(y_1, y_2) \in A_x^2$ such that either $1 \leq y_2 \leq y_1 \leq x$ or $x < y_1 \leq y_2 \leq n$ are counted in the numerator of the summand, we get that this is at least
\begin{align*}
    \frac{\alpha}{2N^2}\sum_{x\in [n],\ A_x\neq\varnothing}\frac{1}{|A_x|}\Biggl(\frac{|A_x \cap [1,x]|^2}{2}+\frac{|A_x \cap [x+1,n]|^2}{2}\Biggr) &\geq \frac{\alpha}{2N^2}\sum_{x\in [n],\ A_x\neq\varnothing}\frac{1}{|A_x|}\cdot \frac{|A_x|^2}{4}\\
    &= \frac{\alpha}{8N^2}\sum_{x\in [n],\ A_x\neq\varnothing}|A_x|\\
    &= \frac{\alpha^2}{32},
\end{align*}
where we made use of convexity in the first inequality. Thus, by combining these calculations with (\ref{fourier representation}), (\ref{main decomposition}) and keeping in mind that $\reallywidehat{S(f_A)}(1) = 0$, we arrive at the desired conclusion. $\qed$

The following auxiliary lemma is a simple consequence of Parseval's identity, but is vital to the success of the $\mathbb{L}^2$ density increment strategy.

\begin{lemma}
\label{parseval bound}
For a subset $A \subseteq [n]^2$, we have the following bound:
$$\sum_{\gamma\in\widehat{G}\setminus\{1\}}\mathbb{E}_{x\in G}|\reallywidehat{\widetilde{1_A}(x,\cdot)}(\gamma)||\reallywidehat{1_A(x,\cdot)}(\gamma)| \leq 1.$$
\end{lemma}
\noindent\textit{Proof.} If $A_x \neq \varnothing$, then Parseval's identity implies that
$$\sum_{\gamma\in\widehat{G}}|\reallywidehat{\widetilde{1_A}(x,\cdot)}(\gamma)||\reallywidehat{1_A(x,\cdot)}(\gamma)| = \frac{1}{S(1_A)(x)}\sum_{\gamma\in\widehat{G}}|\reallywidehat{1_A(x,\cdot)}(\gamma)|^2 = \frac{1}{S(1_A)(x)}\mathbb{E}_{y\in G}|1_A(x,y)|^2 = 1,$$
whereas if $A_x = \varnothing$, the sum trivially vanishes. The conclusion now follows by taking the average over $x\in G$ and interchanging the order of summation. $\qed$

Assuming the second alternative holds in Proposition \ref{main dichotomy}, we could at this point use Lemma \ref{parseval bound} to obtain the existence of a character $\gamma \in \widehat{G} \setminus \{1\}$ such that $|\widehat{S(f_A)}(\gamma)| \geq \frac{\alpha^2}{64}$. We could then carry on basically in the same fashion as in Roth's original proof of his theorem \cite{roth} to obtain a bound of the form $s(n) \ll \frac{n^2}{\log\log n}$. However, as mentioned before, we pursue a different route, which leads to a better bound. 

The following lemma is technical in nature and is a sort of generalisation of Lemma 4 in \cite{green-szemeredi}. It is a device that will be used to convert large Fourier $\mathbb{L}^p$-mass into large Fourier $\mathbb{L}^2$-mass, which in turn serves as input to the density increment argument.

\begin{lemma}
\label{technical lemma}
Let $b_1 \geq b_2 \geq \ldots \geq b_{N-1} \geq 0$ and $\beta > 0$ be such that $\sum_{j=1}^{N-1}b_j^p \geq \beta^p$, $\sum_{j=1}^{N-1}b_j \leq \beta^q$ for some parameters $p > 1$ and $q \in \mathbb{R}$. Then for any $p' \in (1,p)$, there exists a positive integer $m \leq \Bigl\lceil 2^{\frac{1}{p-1}}\beta^{\frac{p(q-1)}{p-1}}\Bigr\rceil$ such that $\sum_{j=1}^{m}b_j \geq cm^{1-\frac{1}{p'}}\beta$, where $c = (2\zeta(p/p'))^{-1/p}$.
\end{lemma}
\noindent\textit{Proof.} To begin, note that for each $k \in [N-1]$ we have
$$kb_k \leq \sum_{j=1}^{k}b_j \leq \sum_{j=1}^{N-1}b_j \leq \beta^q,$$
whence $b_k \leq k^{-1}\beta^q$. Hence, for all $k \in [N-1]$ we get
$$\sum_{j=k+1}^{N-1}b_j^p \leq b_k^{p-1}\sum_{j=k+1}^{N-1}b_j \leq b_k^{p-1}\beta^q \leq k^{1-p}\beta^{pq}.$$
In particular, taking $k = \min\Bigl\{\Bigl\lceil 2^{\frac{1}{p-1}}\beta^{\frac{p(q-1)}{p-1}}\Bigr\rceil, N-1\Bigr\}$, it follows that $\sum_{j=1}^{k}b_j^p \geq \frac{1}{2}\beta^p$. Now if $b_j \leq cj^{-1/p'}\beta$ for all $j \in [k]$, then
$$\sum_{j=1}^{k}b_j^p \leq c^p\beta^p\sum_{j=1}^{k}j^{-p/p'} < c^p\beta^p\zeta(p/p') = \frac{1}{2}\beta^p,$$
which is absurd. Therefore, there must exist $m \in[k]$ such that $b_m \geq cm^{-1/p'}\beta$, whence
$$\sum_{j=1}^{m}b_j \geq mb_m \geq cm^{1-\frac{1}{p'}}\beta,$$
as desired. $\qed$

Following \cite{henriot}, for a set of characters $\Gamma \subseteq \widehat{G}$ and a parameter $\nu \in (0,2]$, we will say that $\Gamma$ is \emph{$\nu$-annihilated} by a set $X \subseteq G$ if $|\gamma(x)-1| \leq \nu$ for all $\gamma \in \Gamma$ and $x \in X$ . The following simple lemma is a variant of Lemma 12 in \cite{henriot}.

\begin{lemma}
\label{annihilation implies large fourier coefficients}
If $\Gamma \subseteq \widehat{G}$ is $\nu$-annihilated by $X$, then $|\widehat{1_X}(\gamma)| \geq \left(1-\frac{\nu^2}{2}\right)\frac{|X|}{N}$ for all $\gamma \in \Gamma$.
\end{lemma}
\noindent\textit{Proof.} If $\gamma(x) = e(x\theta)$ for $\theta \in \mathbb{R}/\mathbb{Z}$, then for all $x \in X$ we have $2|\sin(\pi x\theta)| = |1-e(x\theta)| \leq \nu$, whence
$$\text{Re}\gamma(x) = \cos(2\pi x\theta) = 1-2\sin^2(\pi x\theta) \geq 1-\frac{\nu^2}{2}.$$
Therefore, we have
$$\text{Re}\widehat{1_X}(\gamma) = \mathbb{E}_{x\in G}1_X(x)\text{Re}\gamma(x) \geq \Bigl(1-\frac{\nu^2}{2}\Bigr)\frac{|X|}{N},$$
as required. $\qed$

The following three propositions, whose names are inspired by \S8 of \cite{henriot}, constitute the heart of the proof. The idea is to exploit the second alternative in Proposition \ref{main dichotomy} via the following trichotomy: we either have large Fourier $\mathbb{L}^2$-mass along characters oscillating in the horizontal direction, or large Fourier $\mathbb{L}^2$-mass along characters oscillating in the vertical direction or a large Fourier coefficient at a character oscillating in the vertical direction. These alternatives are then shown to lead to density increments on either an arithmetic progression of rows or an arithmetic progression of columns. 

\begin{proposition}[Horizontal $\mathbb{L}^2$ density increment]
\label{horizontal l2 density increment}
Let $A \subseteq [n]^2$ be a subset of density $\alpha > 0$ such that
\begin{equation*}
    \sum_{\gamma\in\widehat{G}\setminus\{1\}}\Bigl(\mathbb{E}_{x\in G}|\reallywidehat{\widetilde{1_A}(x,\cdot)}(\gamma)||\reallywidehat{1_A(x,\cdot)}(\gamma)|\Bigr)^{\frac{3}{2}} \geq (C\alpha)^{\frac{3}{2}}.
\end{equation*}
Then there exists a non-empty set $\Gamma \subseteq \widehat{G}\setminus\{1\}$ of size $|\Gamma| \ll \alpha^{-3}$ such that, for any set $P \subseteq G$ that $1$-annihilates $\Gamma$,
\begin{equation*}
    \mathbb{E}_{x\in G}\max_{y\in G}\frac{|(A_x-y) \cap P|}{|P|} \gg C|\Gamma|^{1/4}\alpha.
\end{equation*}
\end{proposition}
\noindent\textit{Proof.} Lemma \ref{parseval bound} implies that, for a suitable enumeration $\widehat{G}\setminus\{1\} = \{\gamma_1,\ldots,\gamma_{N-1}\}$, we may apply Lemma \ref{technical lemma} with $b_j = \mathbb{E}_{x\in G}|\reallywidehat{\widetilde{1_A}(x,\cdot)}(\gamma_j)||\reallywidehat{1_A(x,\cdot)}(\gamma_j)|$, $\beta = C\alpha$, $p = \frac{3}{2}$, $q = 0$ and $p' = \frac{4}{3}$ to obtain a set $\Gamma \subseteq \widehat{G}\setminus\{1\}$ of the required size such that
$$\sum_{\gamma \in \Gamma}\mathbb{E}_{x\in G}|\reallywidehat{\widetilde{1_A}(x,\cdot)}(\gamma)||\reallywidehat{1_A(x,\cdot)}(\gamma)| \gg C|\Gamma|^{1/4}\alpha.$$
Hence, if $\Gamma$ is $1$-annihilated by $P$, we obtain using Lemma \ref{annihilation implies large fourier coefficients} that 
\begin{align*}
    \mathbb{E}_{x\in G}\langle \widetilde{1_A}(x,\cdot) \circ 1_P, 1_A(x,\cdot) \circ 1_P\rangle &= \mathbb{E}_{x\in G}\sum_{\gamma\in\widehat{G}}\reallywidehat{\widetilde{1_A}(x,\cdot) \circ 1_P}(\gamma)\overline{\reallywidehat{1_A(x,\cdot) \circ 1_P}(\gamma)}\\
    &= \mathbb{E}_{x\in G}\sum_{\gamma\in\widehat{G}}|\reallywidehat{\widetilde{1_A}(x,\cdot)}(\gamma)||\reallywidehat{1_A(x,\cdot)}(\gamma)||\widehat{1_P}(\gamma)|^2\\
    &\geq \Biggl(\frac{|P|}{2N}\Biggr)^2\sum_{\gamma\in\Gamma}\mathbb{E}_{x\in G}|\reallywidehat{\widetilde{1_A}(x,\cdot)}(\gamma)||\reallywidehat{1_A(x,\cdot)}(\gamma)|\\
    &\gg C\Biggl(\frac{|P|}{N}\Biggr)^2|\Gamma|^{1/4}\alpha.
\end{align*}
On the other hand, this quantity can be upper bounded as follows:
\begin{align*}
    \mathbb{E}_{x\in G}\langle \widetilde{1_A}(x,\cdot) \circ 1_P, 1_A(x,\cdot) \circ 1_P\rangle &\leq \mathbb{E}_{x\in G}\lVert\widetilde{1_A}(x,\cdot) \circ 1_P\rVert_1\lVert 1_A(x,\cdot) \circ 1_P \rVert_{\infty}\\
    &= \mathbb{E}_{x\in G}\lVert\widetilde{1_A}(x,\cdot)\rVert_1 \lVert 1_P\rVert_1\lVert 1_A(x,\cdot) \circ 1_P \rVert_{\infty}\\
    &\leq \frac{|P|}{N}\mathbb{E}_{x\in G}\lVert 1_A(x,\cdot) \circ 1_P \rVert_{\infty},
\end{align*}
whence the desired conclusion follows. $\qed$

\noindent\textbf{Remark.} The reason why we use $\widetilde{1_A}$ instead of $1_A$ lies in the last part of the proof of Proposition \ref{horizontal l2 density increment}, where we upper bound the average additive energy of a vertical slice of $A$ with $P$. Had we used $1_A$ instead of $\widetilde{1_A}$, we would end up with a weighted average having the column densities of $A$ as weights. Since these densities may vary a lot from column to column, it seems to us that this would not imply the desired density increment. 

\begin{proposition}[Vertical $\mathbb{L}^2$ density increment]
\label{vertical l2 density increment}
Let $A \subseteq [n]^2$ be a subset of density $\alpha > 0$ such that
$$\sum_{\gamma\in\widehat{G}\setminus\{1\}}|\widehat{S(f_A)}(\gamma)|^{\frac{5}{2}} \geq (C\alpha)^{\frac{5}{2}}.$$
Then there exists a non-empty set $\Gamma \subseteq \widehat{G}\setminus\{1\}$ of size $|\Gamma| \ll \alpha^{-5}$ such that, for any set $P \subseteq G$ that $1$-annihilates $\Gamma$,
\begin{equation*}
    \max_{y\in G}\frac{1}{|P|}\sum_{x\in P}S(1_A)(x+y) \gg C^2|\Gamma|^{1/6}\alpha.
\end{equation*}
\end{proposition}
\noindent\textit{Proof.} On observing that Parseval's identity implies
$$\sum_{\gamma\in\widehat{G}}|\widehat{S(f_A)}(\gamma)|^2 = \mathbb{E}_{x\in G}|S(f_A)(x)|^2 \leq \lVert S(f_A)\rVert_{\infty}\lVert S(f_A)\rVert_1 \leq \lVert S(1_A)\rVert_1 + \alpha\lVert S(1_{[n]^2})\rVert_1 = \frac{\alpha}{4} + \frac{\alpha}{4} = \frac{\alpha}{2},$$
we may apply Lemma \ref{technical lemma} with $b_j = |\widehat{S(f_A)}(\gamma_j)|^2$ for some $\{\gamma_1, \ldots, \gamma_{N-1}\} = \widehat{G}\setminus\{1\}$, $\beta = (C\alpha)^2$, $p = \frac{5}{4}$, $q = \frac{1}{2}$ and $p' = \frac{6}{5}$ to obtain a set $\Gamma \subseteq \widehat{G}\setminus\{1\}$ of the required size such that
$$\sum_{\gamma\in\Gamma}|\widehat{S(f_A)}(\gamma)|^2 \gg C^2|\Gamma|^{1/6}\alpha^2.$$
To pass from the balanced function to the indicator function, we observe the following inequalities:
\begin{align*}
    \Biggl|\sum_{\gamma\in\Gamma}|\widehat{S(f_A)}(\gamma)|^2-\sum_{\gamma\in\Gamma}|\widehat{S(1_A)}(\gamma)|^2\Biggr| &\leq \sum_{\gamma\in\Gamma}\Bigl||\widehat{S(f_A)}(\gamma)|^2-|\widehat{S(1_A)}(\gamma)|^2\Bigr|\\
    &\leq \sum_{\gamma\in\Gamma}|\widehat{S(f_A)}(\gamma)-\widehat{S(1_A)}(\gamma)|(|\widehat{S(f_A)}(\gamma)|+|\widehat{S(1_A)}(\gamma)|)\\
    &\leq \sum_{\gamma\in\Gamma}\alpha|\reallywidehat{S(1_{[n]^2})}(\gamma)|(\lVert S(f_A)\rVert_1 + \lVert S(1_A)\rVert_1)\\
    &\leq \frac{3}{4}\alpha^2\sum_{\gamma\in\Gamma}|\reallywidehat{S(1_{[n]^2})}(\gamma)|.
\end{align*}
On writing $\Gamma = \{\gamma_1,\ldots,\gamma_m\}$, where $\gamma_j(x) = e(a_jx/N)$ for $j \in [m]$ with $a_1, \ldots, a_m \in [N-1]$, we have the following standard estimate for the Fourier coefficients of $S(1_{[n]^2}) = \frac{1}{2}1_{[n]}$:
$$|\reallywidehat{S(1_{[n]^2})}(\gamma_j)| = \frac{1}{2N}\Biggl|\sum_{x=0}^{n-1}e(a_jx/N)\Biggr| = \frac{1}{2N}\Biggl|\frac{1-e(a_jn/N)}{1-e(a_j/N)}\Biggr| \leq \frac{1}{4\min\{a_j,N-a_j\}}.$$
Hence, we have
$$\Biggl|\sum_{\gamma\in\Gamma}|\widehat{S(f_A)}(\gamma)|^2-\sum_{\gamma\in\Gamma}|\widehat{S(1_A)}(\gamma)|^2\Biggr| \ll \alpha^2 \sum_{j=1}^{\lceil m/2\rceil}\frac{1}{j} \leq \alpha^2(1+\log |\Gamma|),$$
so by choosing $C$ to be large enough, we can ensure that the following holds:
$$\sum_{\gamma\in\Gamma}|\widehat{S(1_A)}(\gamma)|^2 \gg C^2|\Gamma|^{1/6}\alpha^2.$$
We now conclude similarly as in the proof of Proposition \ref{horizontal l2 density increment}. For any $P \subseteq G$ that $1$-annihilates $\Gamma$, we get using Lemma \ref{annihilation implies large fourier coefficients} that
$$\lVert S(1_A) \circ 1_P\rVert_2^2 = \sum_{\gamma\in\widehat{G}}|\widehat{S(1_A)}(\gamma)|^2|\widehat{1_P}(\gamma)|^2 \gg C^2\Biggl(\frac{|P|}{2N}\Biggr)^2|\Gamma|^{1/6}\alpha^2.$$
On the other hand, we have the following upper bound on this quantity:
$$\lVert S(1_A) \circ 1_P\rVert_2^2 \leq \lVert S(1_A)\circ 1_P\rVert_1\lVert S(1_A)\circ 1_P\rVert_{\infty} = \frac{\alpha|P|}{4N}\lVert S(1_A)\circ 1_P\rVert_{\infty}.$$
Thus, the proof is complete. $\qed$

It remains to produce a density increment in the case when there is a non-trivial character $\gamma$ such that $|\widehat{S(f_A)}(\gamma)| \gg \alpha$. This is a standard argument as in the proof of Roth's theorem. We require (the $m = 1$ case of) Dirichlet's approximation theorem (see \cite{schmidt}):

\begin{lemma}[Dirichlet's approximation theorem]
\label{dirichlet's approximation theorem}
Let $\theta_1, \ldots, \theta_m \in \mathbb{R}/\mathbb{Z}$ and let $Q$ be a positive integer. Then there exists $q \in [Q]$ such that $\lVert q\theta_j\rVert_{\mathbb{R}/\mathbb{Z}} \leq \frac{1}{Q^{1/m}}$ for all $j \in [m]$. 
\end{lemma}

Lemma \ref{dirichlet's approximation theorem} also has the following simple corollary, which will be later used to furnish arithmetic progressions $P$ as input to Propositions \ref{horizontal l2 density increment} and \ref{vertical l2 density increment}. In what follows, by the \emph{span} of an arithmetic progression we mean the product of its length and common difference.

\begin{corollary}
\label{progression annihilating a set of characters}
Let $\alpha \in (0,1]$ and let $\Gamma \subseteq \widehat{G}\setminus\{1\}$ be a set of characters of size $m \ll \alpha^{-5}$. Then
\begin{enumerate}[(i)]
    \item either $\alpha \ll (\log n)^{-1/6}$;
    \item or $\Gamma$ is $1$-annihilated by an arithmetic progression $P \subseteq [n]$ of length $|P| \gg \alpha n^{\frac{1}{m+1}}$ and span at most $\alpha n$.
\end{enumerate}
\end{corollary}
\noindent\textit{Proof.} Let $\ell \vcentcolon= \Bigl\lfloor \frac{1}{6}(\alpha n)^{\frac{1}{m+1}}\Bigr\rfloor$ and $Q \vcentcolon= (6\ell)^m$. Here, we are assuming that $\alpha n \geq 6^{m+1}$, which is certainly the case unless (i) holds. Thus, $\ell \geq 1$ and $\ell \geq \frac{1}{12}(\alpha n)^{\frac{1}{m+1}} \gg \alpha n^{\frac{1}{m+1}}$. Enumerating $\Gamma$ as $\{\gamma_1,\ldots,\gamma_m\}$, where $\gamma_j(x) = e(x\theta_j)$ for $j \in [m]$ and $\theta_1, \ldots, \theta_m \in \mathbb{R}/\mathbb{Z}$, Lemma \ref{dirichlet's approximation theorem} implies the existence of $q \in [Q]$ such that $\lVert q\theta_j\rVert_{\mathbb{R}/\mathbb{Z}} \leq \frac{1}{6\ell}$ for all $j \in [m]$. Consequently, if we let $P$ be the arithmetic progression $\ell \cdot [q]$, it follows that $\lVert x\theta_j\rVert_{\mathbb{R}/\mathbb{Z}} \leq \frac{1}{6}$ for all $j \in [m]$ and $x \in P$. Thus, $\Gamma$ is $1$-annihilated by $P$ and $P$ has span $\ell q \leq \alpha n$, as desired. $\qed$

\begin{proposition}[Vertical $\mathbb{L}^{\infty}$ density increment]
\label{vertical linfty density increment}
Let $A \subseteq [n]^2$ be a subset of density $\alpha$. Suppose there exists a non-trivial character $\gamma \in \widehat{G}\setminus\{1\}$ such that $|\widehat{S(f_A)}(\gamma)| \geq 4c'\alpha$, where $c' \in (0,\frac{1}{4})$ is an absolute constant. Then
\begin{enumerate}[(i)]
    \item either $\alpha \ll n^{-1/2}$;
    \item or there is an arithmetic progression $P \subseteq [n]$ of length $|P| \gg \alpha n^{1/2}$ and span at most $2c'\alpha n$ such that
\begin{equation*}
    \frac{1}{n|P|}\sum_{x\in P}|A_x| \geq (1+3c')\alpha.
\end{equation*}
\end{enumerate}
\end{proposition}
\noindent\textit{Proof.} To begin, set $Q \vcentcolon= \lceil 2(\pi n)^{1/2}\rceil \leq 4(\pi n)^{1/2}$ and $\ell \vcentcolon= \lfloor c'\alpha n/Q\rfloor$. Provided $\alpha \geq \frac{4}{c'}(\pi/n)^{1/2}$, which we may assume since otherwise (i) holds, we have $\ell \geq 1$ and hence in particular $\ell \geq \frac{c'\alpha n}{2Q} \gg \alpha n^{1/2}$. Writing $\gamma(x) = e(x\theta)$ with $\theta \in \mathbb{R}/\mathbb{Z}$, Lemma \ref{dirichlet's approximation theorem} implies the existence of $q \in [Q]$ such that $\lVert q\theta\rVert_{\mathbb{R}/\mathbb{Z}} \leq \frac{1}{Q}$. Now partition $[n]$ into residue classes modulo $q$. Since each residue class has size at least $\lfloor n/q\rfloor \geq \ell$, we may partition it into consecutive blocks each of length in $[\ell, 2\ell)$. Letting the sets so obtained be $P_1, \ldots, P_k$, it follows that
$$\frac{1}{N}\sum_{j=1}^{k}\Biggl|\sum_{x\in P_j}S(f_A)(x)\gamma(x)\Biggr| \geq 4c'\alpha.$$
Since $\gamma$ fluctuates by at most $\frac{4\pi\ell}{Q} \leq c'\alpha$ on each $P_j$, it follows that
$$\frac{1}{N}\sum_{j=1}^{k}\Biggl|\sum_{x\in P_j}S(f_A)(x)\Biggr| \geq 3c'\alpha.$$
But since
$$\frac{1}{N}\sum_{j=1}^{k}\sum_{x\in P_j}S(f_A)(x) = 0,$$
we deduce that
$$\frac{1}{n}\sum_{j=1}^{k}\max\Biggl(\sum_{x\in P_j}S(f_A)(x), 0\Biggr) \geq 3c'\alpha.$$
By the pigeonhole principle, there exists $j \in [k]$ such that
$$\frac{1}{|P_j|}\sum_{x\in P_j}S(f_A)(x) \geq \frac{3}{2}c'\alpha.$$
Unfolding what this means in terms of the cardinalities of the vertical slices of $A$, we get that (ii) holds, as desired. $\qed$

Finally, we are ready to assemble all the pieces we have developed so far into a density increment result:

\begin{proposition}
\label{density increment on arithmetic progressions}
Let $A \subseteq [n]^2$ be skew corner-free of density $\alpha$. Then at least one of the following holds:
\begin{enumerate}[(i)]
    \item $\alpha \ll (\log n)^{-1/6}$;
    \item there is an arithmetic progression $P \subseteq [n]$ of length $|P| \gg \alpha n^{1/2}$ and span at most $2c'\alpha n$ such that $\frac{|A \cap (P \times [n])|}{|P\times [n]|} \geq (1+3c')\alpha$, where $c' > 0$ is an absolute constant;
    \item there is a positive integer $m \ll \alpha^{-5}$, an arithmetic progression $P \subseteq [n]$ of length $|P| \gg \alpha n^{\frac{1}{m+1}}$ and span at most $\alpha n$ and a skew corner-free set $A' \subseteq P \times [n]$ of density at least $3m^{1/6}\alpha$;
    \item there is a positive integer $m \ll \alpha^{-3}$, an arithmetic progression $P \subseteq [n]$ of length $|P|\gg \alpha n^{\frac{1}{m+1}}$ and span at most $\alpha n$ and a skew corner-free set $A'' \subseteq [n]\times P$ of density at least $3m^{1/4}\alpha$.
\end{enumerate}
\end{proposition}
\noindent\textit{Proof.} Assume that (i) doesn't hold. In particular, Proposition \ref{main dichotomy} implies that
$$\sum_{\gamma \in \widehat{G}\setminus\{1\}}|\reallywidehat{S(f_A)}(\gamma)|\mathbb{E}_{x\in G}|\reallywidehat{\widetilde{1_A}(x,\cdot)}(\gamma)||\reallywidehat{1_A(x,\cdot)}(\gamma)| \gg \alpha^2.$$
Hence, by applying Hölder's inequality, we obtain that
$$\Biggl(\sum_{\gamma\in\widehat{G}\setminus\{1\}}|\reallywidehat{S(f_A)}(\gamma)|^3\Biggr)^{\frac{1}{3}}\Biggl(\sum_{\gamma\in\widehat{G}\setminus\{1\}}\Bigl(\mathbb{E}_{x\in G}|\reallywidehat{\widetilde{1_A}(x,\cdot)}(\gamma)||\reallywidehat{1_A(x,\cdot)}(\gamma)|\Bigr)^{\frac{3}{2}}\Biggr)^{\frac{2}{3}} \gg \alpha^2.$$
Therefore, we have either
\begin{equation}\label{large vertical sum}
    \sum_{\gamma\in\widehat{G}\setminus\{1\}}|\reallywidehat{S(f_A)}(\gamma)|^3 \gg \alpha^3
\end{equation}
or
\begin{equation}\label{large horizontal sum}
    \sum_{\gamma\in\widehat{G}\setminus\{1\}}\Bigl(\mathbb{E}_{x\in G}|\reallywidehat{\widetilde{1_A}(x,\cdot)}(\gamma)||\reallywidehat{1_A(x,\cdot)}(\gamma)|\Bigr)^{\frac{3}{2}} \geq (C\alpha)^{\frac{3}{2}}.
\end{equation}
If (\ref{large horizontal sum}) holds, we may apply Proposition \ref{horizontal l2 density increment} together with Corollary \ref{progression annihilating a set of characters}. Taking $C$ to be sufficiently large, we obtain an arithmetic progression $P \subseteq [n]$ of length and span as in (iv), and with the property that
\begin{equation*}
    \mathbb{E}_{x\in G}\max_{y\in G}\frac{|(A_x-y) \cap P|}{|P|} \geq \frac{3}{2}m^{1/4}\alpha,
\end{equation*}
where $m \ll \alpha^{-3}$ is a positive integer. But then it is easy to see that we may take $A''$ to be the set $\bigcup_{x\in G}[\{x\} \times ((A_x-y_x) \cap P)]$, where $y_x \in G$ is chosen so as to maximise $\frac{|(A_x-y_x) \cap P|}{|P|}$. 

Now consider the case when (\ref{large vertical sum}) holds. Then either $\sum_{\gamma\in\widehat{G}\setminus\{1\}}|\widehat{S(f_A)}(\gamma)|^{\frac{5}{2}} \geq (C\alpha)^{\frac{5}{2}}$ or there exists $\gamma \in \widehat{G}\setminus\{1\}$ such that $|\widehat{S(f_A)}(\gamma)| \gg \alpha$. In the former case, Proposition \ref{vertical l2 density increment} combined with Corollary \ref{progression annihilating a set of characters} yields an arithmetic progression $P \subseteq [n]$ of length and span as in (iii) such that
$$\max_{y\in G}\frac{1}{|P|}\sum_{x\in P}S(1_A)(x+y) \geq \frac{3}{2}m^{1/6}\alpha,$$
where $m \ll \alpha^{-5}$ is a positive integer. But we may now take $A'$ to be the set $\bigcup_{x\in P}(\{x\}\times A_{x+y})$, where $y \in G$ is chosen so as to maximise $\frac{1}{|P|}\sum_{x\in P}S(1_A)(x+y)$. Finally, in the remaining case, Proposition \ref{vertical linfty density increment} immediately implies that (ii) holds. This concludes the proof. $\qed$

However, we need a density increment on a square substructure. Luckily, this can be achieved by a simple averaging argument:

\begin{lemma}
\label{pigeonholing lemma}
Let $P \subseteq [n]$ be an arithmetic progression of length $\ell$ and common difference $d$. If $B \subseteq P \times [n]$ is a subset of density $\beta$, then there is a translate $P' \subseteq [n]$ of $P$ such that $B$ has density at least $\beta - \frac{\ell d}{n}$ on $P \times P'$. Similarly, if $B \subseteq [n] \times P$ has density $\beta$, then there is a translate $P' \subseteq [n]$ of $P$ such that $B$ has density at least $\beta - \frac{\ell d}{n}$ on $P' \times P$.
\end{lemma}
\noindent\textit{Proof.} We only prove the first statement; the proof of the second statement is completely analogous. Partition $[\ell d]$ into $d$ translates of $P$ and then partition $[n]$ into translates of $[\ell d]$ and an error set $E$ of size less than $\ell d$. Then $P \times E$ occupies at most $\ell^2 d$ elements of $B$, so $B \setminus (P\times E)$ has density at least $\beta - \frac{\ell d}{n}$ on $P \times [n]$. The conclusion now follows by averaging. $\qed$

By putting together Proposition \ref{density increment on arithmetic progressions} and Lemma \ref{pigeonholing lemma}, we obtain the final form of our density increment result:

\begin{proposition}
\label{final density increment}
Let $A \subseteq [n]^2$ be skew corner-free of density $\alpha$. Then
\begin{enumerate}[(i)]
    \item either $\alpha \ll (\log n)^{-1/6}$;
    \item or there exist positive integers $m \ll \alpha^{-5}$, $n' < n$ such that $n' \gg \alpha n^{\frac{1}{m+1}}$ and there exists a skew corner-free set $A' \subseteq [n']^2$ of density at least $(1+c')m^{1/6}\alpha$, where $c' > 0$ is an absolute constant.
\end{enumerate}
\end{proposition}

To establish Theorem \ref{main upper bound}, all that remains is to iterate Proposition \ref{final density increment}. We will show that the bound in fact holds with $c_2 = \min(c'/3,1/6)$. To this end, we show by induction on $n$ that any skew corner-free subset of $[n]^2$ for $n > 1$ has density at most $C'(\log n)^{-c_2}$, where $C' > 0$ is some large constant. In case (i) holds in Proposition \ref{final density increment}, we are done by taking $C'$ to be sufficiently large. Otherwise, there exist positive integers $m, n'$ as in (ii) and a skew corner-free subset of $[n']^2$ of density at least $(1+c')m^{1/6}\alpha$. Here, we may assume that $n$ is sufficiently large, that $m \leq \frac{1}{24}(\log n)^{5/6}$ and $n' \geq n^{\frac{1}{m+1}}(\log n)^{-1/6} > 1$. Hence, by the induction hypothesis, we have $(1+c')m^{1/6}\alpha \leq C'(\log n')^{-c_2}$. Assuming for the sake of contradiction that $\alpha > C'(\log n)^{-c_2}$, it follows that
\begin{equation}\label{final contradiction}
    (1+c')^{\frac{1}{c_2}}m^{\frac{1}{6c_2}}\frac{\log n'}{\log n} < 1.
\end{equation}
But $\frac{\log n'}{\log n} \geq \frac{1}{m+1} - \frac{\log\log n}{6\log n}$ and $\log n \geq (24m)^{6/5}$. Since the function $x \mapsto \frac{\log x}{x}$ is decreasing on the interval $[e,\infty)$, it follows that $\frac{\log\log n}{6\log n} \leq \frac{\frac{1}{5}\log(24m)}{(24m)^{6/5}}$, which in turn is at most $\frac{1}{4m}$. Indeed, this follows from the inequality $30m^{1/5} \geq \log(24m)$, which holds for all $m \geq 1$. Thus, $\frac{\log n'}{\log n} \geq \frac{1}{4m}$, so the left-hand side in (\ref{final contradiction}) is at least
$$\frac{1}{4}\Bigl(1+\frac{c'}{c_2}\Bigr)m^{\frac{1}{6c_2}-1} \geq 1,$$
which is a contradiction. This concludes the proof.

\section{Concluding remarks and open problems}

The proofs of both Theorem \ref{main lower bound} and Theorem \ref{main upper bound} exhibit many parallels between the problem of estimating $s(n)$ and that of determining the largest size of a subset of $[n]$ lacking three-term arithmetic progressions. The latter problem is one of the cornerstones of additive combinatorics and has been a subject of intense study in the last seventy years. Recently, a huge breakthrough on this problem was made by Kelley and Meka \cite{kelley-meka}, who proved upper bounds matching the shape of the well-known lower bound of Behrend \cite{behrend}. This came as somewhat of a surprise since Kelley and Meka largely abandon the Fourier-analytic approach, relying instead on methods based in physical space. Thus, it is an interesting challenge to investigate the scope of these methods and explore potential applications to other related problems. In the very recent survey \cite{peluse-finite-fields}, some optimism has been expressed as to the possibility of applying Kelley-Meka techniques to improve bounds for corner-free sets. Given the apparently even larger similarity of skew corners to three-term arithmetic progressions, we think that a good intermediate step would be to investigate whether these techniques could be used to the improve bounds for the problem at hand. Specifically, we ask the following question:

\begin{question}
\label{kelley-meka bounds for skew corners}
Is it true that $s(n) \ll \frac{n^2}{\exp(c(\log n)^{\beta})}$ for some constants $c, \beta > 0$?
\end{question}

There is another rather intriguing parallel with the study of corner-free sets which we think deserves more attention. As mentioned before, the usual notion of Fourier uniformity is not sufficient to get a handle on the count of corners in a given set. A concrete example explaining this phenomenon is given by Example 5.3 in \cite{green-finite-fields}, which also appears in \S1 of \cite{shkredov}. In this example, a random set of $B \subseteq \mathbb{F}_2^n$ of density $\beta$ is taken and the product set $A = B^2 \subseteq (\mathbb{F}_2^n)^2$ is considered. Being the product of a random set with itself, $A$ is a highly uniform set of density $\alpha = \beta^2$. However, the count of corners in $A$ is equal to the number of triples $(x,y,d)\in(\mathbb{F}_2^n)^3$ such that $x,y,x+d,y+d \in B$, which is $\beta^4N^3 = \alpha^2N^3$ in expectation. This is significantly more than $\alpha^3N^3$ corners, what one would expect for a truly random subset of $(\mathbb{F}_2^n)^2$ of density $\alpha$.

Taking a closer look at this example, it becomes clear that a similar reasoning applies to our problem: if $B \subseteq \mathbb{F}_2^n$ is random of density $\beta$, then $A = B^2$ contains roughly $\beta^5N^4 = \alpha^{5/2}N^4$ skew corners, whereas a truly random subset of density $\alpha$ would contain more like $\alpha^3N^4$ skew corners. However, this apparent deficiency did not prevent us from using Fourier analysis to obtain bounds for skew corner-free sets. Indeed, we have shown in Proposition \ref{generalised von neumann} that any uniform subset has roughly \emph{at least} as many skew corners as a random set of the same density. In light of this, we do not find the example discussed in the previous paragraph as entirely satisfactory evidence as to the impossibility of applying Fourier analysis to the study corner-free sets. This motivates us to ask the following question:

\begin{question}
\label{uniform set with fewer than expected corners}
Do there exist Fourier uniform sets with significantly fewer than expected corners? More precisely, do there exist $\alpha, \delta > 0$ such that, for infinitely many $N$, there exists a set $A \subseteq (\mathbb{Z}/N\mathbb{Z})^2$ of density $\alpha+o(1)$ which is $o(1)$-uniform, but has at most $(\alpha^3 - \delta)N^3$ corners?
\end{question}

We believe that Question \ref{uniform set with fewer than expected corners} should have an affirmative answer. However, showing the existence of such sets might not be a simple task. Indeed, comparing this to the situation for four-term arithmetic progressions, it has been known since the seminal work of Gowers \cite{gowers-szemeredi} that there exist uniform sets with \emph{more} than expected arithmetic progressions of length four. In spite of this, constructing such sets with \emph{fewer} than expected $4$-progressions turned out to be a non-trivial problem on its own (see \cite{gowers-uniform}).

\bigskip

\noindent\textbf{Acknowledgements.} This work was supported by the Croatian Science Foundation under the project number HRZZ-IP-2022-10-5116 (FANAP). The author is grateful to Vjekoslav Kovač and Rudi Mrazović for bringing \cite{pratt} to his attention, as well as for fruitful discussions and encouragement. He would also like to thank Cosmin Pohoata and Dmitrii Zakharov for helpful comments on a draft of the paper.

\bibliographystyle{plain}
\bibliography{references}

\end{document}